\def\Q{{\mathbb Q}}
\def\Z{{\mathbb Z}}

\def\B{{\cal B}}
\def\W{{\cal F}}
\def\F{\hbox{\bf F}}
\def\Fq{\F_q}
\def\mod#1{\,({\rm mod\ }#1)}
\def\wedge{\triangleright}

\def\proof{\noindent {\it Proof.}\hskip 10pt}
\def\endproof{\hfill\vbox{\hrule
    \hbox{\vrule\kern4pt\vbox{\kern4pt
    \kern4pt}\kern4pt\vrule}\hrule}\bigskip}
\documentclass[11pt]{article}
\usepackage{epsfig}
\usepackage{amssymb}
\newtheorem{lemma}{Lemma}[section]

\newtheorem{proposition}[lemma]{Proposition}
\newtheorem{theorem}[lemma]{Theorem}

\begin{document}
\title{Self-interacting polynomials}
\author{Franco Vivaldi}
\date{\it\small
School of Mathematical Sciences, Queen Mary, University of London, \\
London E1 4NS, UK
}
\maketitle
\begin{abstract}
We introduce a class of dynamical systems of algebraic origin,
consisting of self-interacting irreducible polynomials over a field.
A polynomial $f$ is made to act on a polynomial $g$ by mapping the 
roots of $g$. This action identifies a new polynomial $h$, as 
the minimal polynomial of the displaced roots. 
By allowing several polynomials to act on one another, we obtain 
a self-interacting system with a rich dynamics, which affords a 
fresh viewpoint on some algebraic dynamical constructs. 
We identify the basic invariant sets, and study in some detail
the case of quadratic polynomials. We perform some experiments
on self-interacting polynomials over finite fields.

\end{abstract}

\section{Introduction}

From the viewpoint of algebraic dynamics, a polynomial with 
coefficients in a field $K$ can be interpreted in two ways.
On the one hand, it defines a dynamical system over the 
algebraic closure $\bar K$ of $K$; on the other, its roots 
define elements of $\bar K$. 
In the former case the polynomial represents a function acting 
on a state space, that is, a dynamical system; in the latter, 
the polynomial encodes points of the state space itself.

In this paper we exploit this duality to define a new 
class of discrete dynamical systems resulting from 
the mutual interaction between irreducible polynomials.
Such interaction features an active element 
(a function) and a passive one (a set of points, which are
algebraic conjugates), whose role is interchangeable. 
Specifically, given two irreducible polynomials $f$ and $g$ 
over the same field, we define $h=f\wedge g$ to be the minimal 
polynomial of the image of the roots of $g$ under $f$.
This relationship between $f,g$ and $h$ defines the 
{\it wedge operator\/} $\wedge$ (see section \ref{section:SelfInteraction} for details). 
The action of the wedge operator is quite natural, since it performs 
simultaneously the time-evolution of all algebraic conjugates, 
which, algebraically speaking, are indistinguishable.  
This construction, introduced in \cite{Vivaldi}, has been 
developed further in
\cite{BatraMorton,BatraMortonII,CohenHachenberger:99,CohenHachenberger:00}.

The idea is to start from a collection of polynomials,
and then allow them to act on one another via the wedge
operator, thereby producing a new set of polynomials. 
Unlike in conventional constructions, here
the state space generates its own dynamics. 
Our interest in such systems is motivated by the desire 
of developing an algebraic variant of some abstract models 
of adaptive systems and chemical interactions, where functions 
were made to act on other functions via composition, within
the framework of $\lambda$-calculus
\cite{Fontana:90,FontanaBuss:94}.
Subsequently, a similar concept was proposed
as functional dynamics on coupled map lattices
\cite{KataokaKaneko:00,KataokaKaneko:01} using
again the composition of smooth functions as interaction.

The dynamics of self-interaction presents considerable 
difficulties: its phenomenology is overwhelming, while 
the theory, still in an embryonic stage, has yet to 
develop adequate links with other areas of dynamics.
The purpose of this paper is to explore some aspects of 
this problem in an algebraic context. 
Thus the basic question of invariance under 
self-interaction will lead to periodic and pre-periodic 
orbits polynomials, their discriminants and Galois groups, 
while the simplest instances of periodicity brings 
connections with periodicity in some strongly chaotic systems 
(the iterated monomial maps, which lead to roots of unity).

Specifically, we now define the {\it self-image\/} of a set of 
irreducible polynomials to be the result of all mutual nontrivial
wedge interactions, that is,
\begin{equation}\label{eq:SelfImage}
\W:\{f_i\}\,\mapsto\,\left\{f_i\wedge f_j\right\}\hskip 50pt
i,j\in J\qquad i\not=j
\end{equation}
where $J$ is a set of indices. 
We have excluded the trivial action of a polynomial on itself 
---see below--- while the choice of global coupling (as opposed 
to, say, nearest neighbours), eliminates topological considerations 
and is justified by the fact that we will be dealing mostly with 
small sets.

\begin{figure}[h]
\hfil\epsfig{file=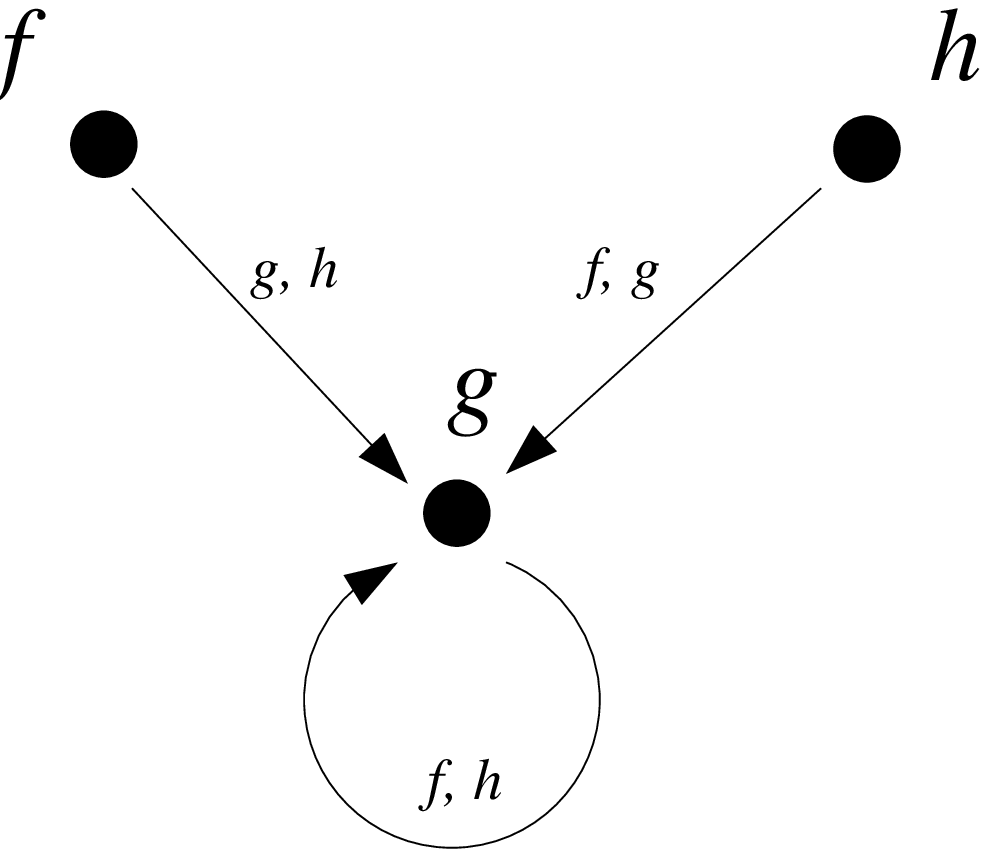,width=5cm,height=4cm}\hfil
\caption{\label{fig:Stable3Set} 
The graph of a stable 3-set $\{f,g,h\}$ over the finite 
field $\F_3$, where $f(x)=x^2+1,\,\,\, g(x)=x^2+x+2,\,\,\, h(x)=x^2+2x+2$.
The notation $g\stackrel{f}{\rightarrow}h$ indicates that
$f$ maps the roots of $g$ into those of $h$; accordingly, 
this graph represents the six equations
$h\wedge g=g,\, f\wedge g=g,\, f\wedge h=g,\, h\wedge f=g,\,
g\wedge h=g,\, g\wedge f=g.$
}
\end{figure}

As a phase space for the dynamical system generated by $\W$, 
we can take any superset of $\{f_i\}$ with the property of
being {\it stable,} that is, closed with respect to the wedge operator.
It is expedient to arrange the polynomials of a stable set
as vertices of a {\it polynomial graph\/} $\Gamma$
\cite{BatraMorton} (see figure \ref{fig:Stable3Set}).
Two polynomials $g$ and $h$ are joined by an oriented arc
$(g,h)$ if there exists a polynomial $f$ among the 
vertices sending the roots of $g$ onto the roots of 
$h$, that is, if $f\wedge g=h$.

In this work we confine our attention to a restricted range of 
problems: classifying small stable sets, analyzing simple
instances of periodic behaviour, and identifying algebraic 
mechanisms leading to highly organized self-images.
In the next section we define precisely the wedge operator,
stable sets and their graphs, and classify stable 2-sets in 
terms of right and left invariants of the wedge operator
---see figure \ref{fig:Stable2Sets}.
We find that stable 2-sets originate from the existence of
certain relations between a polynomial, its periodic points
or its pre-images of zero (or both).
In section \ref{section:Blocks} we partition polynomials 
into {\it blocks}, which are obtained by shifting the roots 
of a polynomial by the elements of the ground field. 
Blocks turn out to have organized self-images, and the 
self-intersection of blocks leads to appearance of 
right-left invariants of the wedge operator 
(theorem \ref{theorem:RightLeftInvariants}).

Quadratic polynomials are studied in detail in
section \ref{section:Quadratic}, resulting in
a fairly complete characterization of their
graphs (theorem \ref{theorem:Quadratic}), and
of the dynamics of quadratic 2-sets, in terms of a 
three-dimensional skew map over the ground field.
We show that periodic 2-sets originate from certain
roots of unity, in a field-theoretic setting
not dissimilar to that of monomial maps
(section \ref{subsection:Periodic2Sets}).

Some aspects of self-interaction over finite fields 
are dealt with in section \ref{section:FiniteFields}.
We count the number of blocks (theorem \ref{theorem:NumberOfBlocks}), 
and, for quadratic polynomials and odd characteristics, the 
number of stable and periodic 2-sets
(theorems \ref{theorem:NumberOfStable2Sets} and 
\ref{theorem:NumberOfPeriodic2Sets}).
We conclude by describing the results of some experiments
over finite fields, concerning the occurrence of
small stable sets, and one asymptotic problem on the 
periodicity of quadratic 2-sets.

\bigskip
I am much indebted with Patrick Morton, for many illuminating 
discussions on various matters concerning this work. 

\section{The wedge equation}
\label{section:SelfInteraction}

Let $K$ be a field, let $f$ and $g$ be monic (unit leading 
coefficient) irreducible polynomials over $K$.
Let $\alpha$ be a root of $g$ and let $K(\alpha)$ be
the field obtained by adjoining $\alpha$ to $K$.
We assume that $K(\alpha)$ is separable (the roots of $g$ 
are distinct).

We represent the minimal polynomial $h$ of $f(\alpha)$ as $h=f\wedge g$; 
equivalently, irreducible $h=f\wedge g$ if and only if $g(x)$ 
divides $h(f(x))$. 
We will show (see below and \cite{Vivaldi})
that $f\wedge g$ is a monic polynomial over $K$, which is 
uniquely defined by $f$ and $g$ ({\it i.e.}, it is independent 
of the choice of the root $\alpha$ of $g$), and whose degree 
is a divisor of the degree of $g$.
Specifically,
\begin{equation}\label{eq:Wedge}
g(x)=\prod_\alpha\,(x-\alpha)
\hskip 30pt
\Longrightarrow
\hskip 30pt
(f \wedge g)(x)=
  \mathop{{\prod}'}_{\alpha\hskip 1pt}\,\left(x-f(\alpha)\right)
\end{equation}
where the primed product is taken over any maximal set 
of roots $\alpha$ of $g$ for which all terms are distinct.
As $(f \wedge g)(0)$ is the resultant of $f$ and $g$
\cite[chapter 1.4]{LidlNiederreiter},
the coefficients of $f \wedge g$ may be regarded as
generalizations of the resultant.

\subsection{Construction of $f\wedge g$}
\label{subsection:Wedge}

We now describe an algorithm to construct the polynomial 
$h=f\wedge g$ from $f$ and $g$ without computing roots, 
using linear algebra (see {\cite{Vivaldi} for details).
Let $g$ have degree $n$, and let $\alpha$ be one of its roots. 
The extension $K(\alpha):K$ has degree $n$, since $g$ is irreducible.
The field $K(\alpha)$ is represented canonically as the quotient 
of the polynomial ring $K[x]$ modulo the maximal ideal 
$(g)$ \cite[page 32]{Waerden}
\begin{equation} \label{eq:Canonical}
K(\alpha)\cong {K[x]/(g(x))}, 
\end{equation}
meaning that the polynomials in $K[x]$ are replaced by their 
remainder upon division by $g(x)$. 
Under this isomorphism, the root $\alpha$ is identified with
the residue class of $x$.
The polynomials $1,x,x^2,\ldots,x^{n-1}$ form a basis of 
$K[x]/(g(x))$ over $K\,(\cong K[x]/(x))$. 
We consider the map of $K[x]/(g(x))$ into itself determined 
by multiplication by $f(x)$. 
This mapping is linear \cite[chapter 2]{Marcus}, and in the above 
basis it is represented by the matrix $M=\{m_{i,j}\}$ over $K$ 
given by
\begin{equation}\label{eq:M}
f(x)\cdot x^{k-1}\,\equiv\, m_{1,k}+m_{2,k}x+\cdots +m_{n,k}x^{n-1}
                           \mod{g(x)}\hskip 30pt k=1,\ldots,n.
\end{equation}
It can be shown that $h$ is the minimal polynomial of $M$, and
that the Jordan form of $M$ has $r$ identical blocks along the 
diagonal, for some integer $r$.
This means that the characteristic polynomial of $M$ takes
the form $\left(h\right)^r$, and hence the degree of $h$ is 
a divisor of that of $g$.
Such divisor is proper whenever the action of $f$ on the 
roots of $g$ is not injective, or, equivalently, whenever 
$K(f(\alpha))$ is a proper subfield of $K(\alpha)$, for 
each root $\alpha$ of $g$.

\subsection{Stable sets and their graphs}
\label{subsection:StableSets}

A set $S$ of irreducible polynomials is said to be 
{\it stable\/} if it contains its self-image $\W(S)$,
defined via the wedge operator in (\ref{eq:SelfImage}). 
Because a polynomial over a field $K$ leaves every 
algebraic extension of $K$ invariant, ``large'' stable 
sets are naturally constructed by means of field extensions.
In particular, the set of minimal polynomials of all
the elements of a finite extension $K(\alpha):K$ is stable;
we denote it by $S(f)$, where $f$ is any irreducible monic
polynomial with a root  $\beta$ such that $K(\beta)=K(\alpha)$.
One sees that this definition does not distinguish between 
conjugate fields.
A prominent subset of $S(f)$ is 
\begin{equation}\label{eq:E}
E(f):=\{g\in S(f)\,|\, S(g)=S(f)\}
\end{equation}
which contains only the polynomials generating the top field.
The set $E(f)$ is not stable, in general.
Furthermore, given $f,g\in K[x]$, the sets $E(f)$ and $E(g)$ 
are either disjoint or one is contained in the other.
An example is given in figure \ref{fig:Stable3Set},
where we display all quadratic irreducible polynomials
over $\F_3$; in this case $E(f)=E(g)=E(h)=\{f,g,h\}$.

The {\it graph\/} $\Gamma(S)$ of a stable set $S$ is the 
directed graph whose vertices are the elements of $S$, and 
where there exists an arc $(g,h)$ from 
$g$ to $h$ if $f\wedge g=h$ for some $f$ in $S$.
For each arc $(g,h)$ we consider the collection $\omega$
of all polynomials $f$ such as $f\wedge g = h$.
The cardinality of $\omega$ will be called the {\it multiplicity\/} 
of the arc $(f,g)$.
To make $\omega$ explicit, we write 
$g\stackrel{\omega}{\rightarrow} h$ in place of $(g,h)$,
see figures \ref{fig:Stable3Set}--\ref{fig:Invariant3Set}.

Graphs arising from polynomials of degree one are straightforward:
for all $a,b\in K$, the polynomial $f(x)=x+b-a$ is the unique 
polynomial sending the root of $g(x)=x-a$ to that of $h(x)=x-b$.
Thus $\Gamma$ is a complete graph, and all arcs have multiplicity one.
The quadratic case is already non-trivial 
---see section \ref{section:Quadratic}. 
In general case, some information on $\Gamma(S)$ is obtained by 
restricting vertices to the subsets $E(f)$ of $S$ defined in 
(\ref{eq:E}); the resulting subgraphs will be called {\it extension graphs.}
In an extension graph of degree $n$, the multiplicity of each arc
is at most $n$. This is a consequence of the following finiteness result

\begin{proposition}\label{proposition:Spectrum}
Let $g$ and $h$ be monic irreducible polynomials 
of positive degree $n$ and $m$, respectively.
Then, for every integer $d$, with $1 \leq d \leq n$, 
there are at most $m$ polynomials $f$ of degree $d$,
such that $f\wedge g=h$.
\end{proposition}

\proof
Let $g,f_1,f_2 \in K[x]$ be non-constant monic polynomials, with 
$g$ irreducible and ${\rm deg}\,f_2 \leq {\rm deg}\, f_1\leq n$. 
Let $\alpha$ be a root of $g$, and let $f_1(\alpha)=f_2(\alpha)$.
We begin to show that if $f_1$ and $f_2$ have the same degree, 
they coincide; otherwise, ${\rm deg}\,f_1=n$, and $g=f_1-f_2$.

Let $l=f_1-f_2$. 
Then $l(\alpha)=0$, and since ${\rm deg}\,l\leq n$ 
and $g$ is the minimal polynomial of $\alpha$, we have 
either $l=0$ or ${\rm deg}\,l=n$, with $l$ a constant multiple of $g$.
If $f_1$ and $f_2$ have the same degree, then the degree of $l$ is lower 
(because $f_1$ and $f_2$ are monic), so that $l=0$, and $f_1=f_2$.
If their degrees are different, then $f_1$ must have degree $n$, 
and since $f_1$ is monic, we have $l=g$ as claimed.

Let now $\beta=f(\alpha)$ be a root of $h$, with
$f$ a monic polynomial of degree $d\leq n$.
We have shown that such $f$, if it exists, is unique.
There are $m$ possible choices of $\beta$ among the
roots of $h$, giving as many polynomials, not 
necessarily irreducible.
The transitivity of the Galois group of $g$ means that
the action of $f$ on the roots of $g$ is uniquely 
determined by the image of the root $\alpha$, so there
are no other polynomials. 
This completes the proof.
\endproof 

\subsection{Wedge invariants and stable 2-sets}\label{section:Invariance}

The search for small stable sets begins from the study of 
right and left invariance under the wedge operator
\begin{equation} \label{eq:Invariants}
 (i)\quad f\wedge g=g\hskip 50pt 
 (ii)\quad g\,\wedge\, f =g.
\end{equation}
In equations $(i)$ and $(ii)$, we say that $g$ is a 
{\it right} and a {\it left wedge invariant} of $f$, 
respectively. These invariants turn out to have dynamical 
significance. 

{\it Right invariants.} \/ They originate from periodic points.
To see this we note that if $g$ is a right invariant 
of $f$, then $f$ maps the roots of $g$ into themselves. 
In fact from the irreducibility of $g$ and the transitivity 
of its Galois group, it is easy to see that $f$ must permute 
the roots of $g$, that is, the roots of $g$ are periodic orbits of $f$
(for Galois theory of periodic orbits see 
\cite{VivaldiHatjispyros,MortonPatel}).
To construct solutions $g$ of (\ref{eq:Invariants}$i$),
we let $f^n$ be the $n$-th iterate of $f$ (with $f^0(x)=x$).
Then the roots of the polynomial $f^n(x)-x$ are the periodic 
points of $f$ of period dividing $n$; to eliminate spurious periods,
we consider the polynomials
\cite{VivaldiHatjispyros,Bousch,MortonPatel}
\begin{equation} \label{eq:PhiPoly}
\Phi_{n,f}(x)\, =\, \prod_{d\,|n}(f^d(x)-x)^{\mu(n/d)}\qquad\quad n=1,2,\ldots
\end{equation}
where $\mu$ is the M\"obius function \cite[chapter 2]{Apostol}. 
The roots of $\Phi_{n,f}$ are periodic points of $f$ of 
{\it essential\/} period $n$ \cite{MortonSilverman},
which is the minimal period apart, possibly, at bifurcations
\cite{MortonVivaldi}.
It follows that every right invariant of $f$ is an irreducible 
factor of $\Phi_{n,f}$, for some $n$. The converse is not true:
if $g$ is an irreducible factor of $\Phi_{n,f}$, all we can say 
is that $g$ is a right invariant of $f^d$, for some divisor
$d$ of $n$.
In particular, if $\Phi_{n,f}$ is irreducible, then 
$f\wedge\Phi_{n,f}=\Phi_{n,f}$. 
(The polynomial $\Phi_{n,f}$ is `generically' irreducible,
that is, irreducible over $K=\bar \Q$, when the coefficients 
of $f$ are regarded as indeterminates \cite{Morton:96}.)

{\it Left invariants.} \/
They originate from pre-images of zero.
Let us consider the functional equation 
$f^k\wedge f^n=f^{n-k},\quad n\geq k$,
ignoring irreducibility for the moment
(the polynomial $f^n$ is generically irreducible over 
any field of characteristic zero \cite{Odoni:85}).
Letting $n=2k$ and $g=f^k$, we obtain $g\wedge g^2=g$. 
Thus, $g$ is a left invariant of any irreducible factor
of $g^2$. Conversely, if $g$ is a left invariant of $f$,
and $\alpha$ is a root of $f$, then $g(g(\alpha))=0$,
that is, $f$ is an irreducible factor of $g^2$.
This construction is straightforward in the degree one case:
for all $0\not=a\in K$, the polynomial $x+a$ is a left 
invariant of $x+2a$.

Right and left invariants are not completely unrelated; 
in section \ref{section:Blocks} we shall establish a 
sufficient condition for their simultaneous occurrence 
(theorem \ref{theorem:RightLeftInvariants}).

Three different types of stable 2-sets can be constructed 
by combining right and left invariants in all possible ways, 
as illustrated in figure \ref{fig:Stable2Sets}. 
Type I and II stable 2-sets consist of two right and two 
left invariants, respectively; they are fixed points of 
the mapping $\W$ --- see equation \ref{eq:SelfImage}. 
The case of mixed invariants is denoted as type III; its
self-image is the right-left invariant.

\begin{figure}
\vskip -50 pt
\hfil\epsfig{file=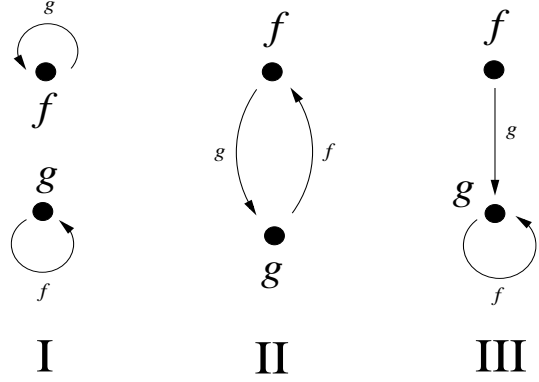,width=7cm,height=5cm}\hfil
\caption{\label{fig:Stable2Sets} 
Stable 2-sets $\{f,g\}$. 
Type I: two right invariants ($f\wedge g=g$ and $g\wedge f=f$). 
Type II: \/ two left invariants ($f\wedge g=f$ and $g\wedge f=g$).
Type III: \/ mixed invariants ($f\wedge g=g$ and $g\wedge f=g$).
Types I and II sets are invariant under $\W$,
while type III collapses onto $\{g\}$.
}
\end{figure}

We will arrive to a complete description of quadratic stable 2-sets in
sections \ref{section:Quadratic} and \ref{section:FiniteFields}.
Type I and III are easily constructed, while type II do not exist,
requiring at least cubic polynomials 
(theorem \ref{theorem:NoQuadraticStableSetsOfTypeII}).

A specific family of type I stable sets can be constructed, involving
polynomials of arbitrarily large degree. Consider a non-constant monic polynomial 
$f$ over $K$. Then the sets
\begin{eqnarray}
\{f(x^2),\,f(x^2)-x\} && \hbox{\rm char}(K)\not=2 \label{eq:TypeI}\\
\{f(x),\,f(x)-x\} && \hbox{\rm char}(K)=2\label{eq:TypeI2}
\end{eqnarray}
are stable of type I whenever the polynomials involved are 
irreducible (char$(K)$ denotes the characteristic of 
$K$, see \cite[chapter 1.2]{LidlNiederreiter}).

To see this, let $g(x)=f(x^2)$, with char$(K)\not=2$.
Then $g(x)-x=\Phi_{1,g}(x)$ is a right invariant of $g$.
Furthermore, if $\alpha$ is a root of $g$, then 
$\Phi_{1,g}(\alpha)=-\alpha$ is a root of $\Phi_{1,g}\wedge g$. 
Therefore, for $g$ to be a right invariant of $\Phi_1$, it suffices 
that $-\alpha$ also be a root of $g$, which is true since $g(x)=f(x^2)$.
It follows that the set (\ref{eq:TypeI}) is a type I stable set 
(pending irreducibility, that is);
the first polynomial induces the identity permutation of the root
of the second, while the second induces 2-cycles on the roots of the first.
In particular, $g$ is an irreducible factor of $\Phi_{2,h}$, where $h=\Phi_{1,g}$. 
If char$(K)=2$, the polynomials in (\ref{eq:TypeI2}) are 
the $\Phi_1$-polynomial of each other.

For example, the polynomial $f(x)=x^2+1$ gives rise to the following 
stable 2-set of type I over $\Q$
$$
\{x^4+1,\,x^4-x+1\}.
$$

\section{Blocks}
\label{section:Blocks}

A {\it block\/} is a maximal collection of polynomials whose 
roots differ by elements of the ground field. 
Blocks turn out to have ``small'' and well-structured self-images. 
These will be studied in section \ref{section:SelfImageOfBlocks},
where we derive a sufficient condition for degree invariance of the 
self-image of a block (proposition \ref{proposition:DegreeInvariance}),
and we show that intersection of a block with its self-image leads 
to right and left invariants of the wedge operator 
(theorem \ref{theorem:RightLeftInvariants}).
Before that, we establish the notation and provide some preliminary 
results. Since the polynomials in a block share the same discriminant, 
the study of discriminants will become relevant.

\subsection{Preliminaries}
\label{section:Preliminaries}
We begin with the group $G$ of matrices
$$
\sigma_{a,b}\,=\,\left(\matrix{a & b \cr 0 & 1 \cr}\right)\qquad 
 a,b \in K,\quad a\not= 0
$$
which we let act on irreducible polynomials as follows
\begin{equation}\label{eq:ActionOfG}
\sigma(f)(x)=a^{-{\rm deg\/}f} f(a\,x+b).
\end{equation}
The polynomial $\sigma(f)$ is monic and irreducible, and since
$\sigma_{c,d}(\sigma_{a,b}(f))=(\sigma_{c,d}\sigma_{a,b})(f)$,
equation (\ref{eq:ActionOfG}) defines an action of $G$ on 
the set $E(f)$, which was defined in (\ref{eq:E}).
The orbits of $G$ partition the vertices of 
every extension graph into {\it clusters.} 

We consider the additive subgroup 
\begin{equation} \label{eq:SigmaPlus}
G^+=\{\sigma_{1,b}\}
\hskip 40pt
\sigma_{1,b}\,=\,\left(\matrix{1 & b \cr 0 & 1 \cr}\right)\qquad b\in K.
\end{equation}
Its orbits are {\it blocks,} which subdivide each cluster.
We denote blocks by $\Theta$ and the block containing $f$
by $\Theta_f$.
Introducing the short-hand notation
$$
f^+_b:=\sigma_{1,b}(f)
$$
we have $\Theta_f\,=\,\{f_b^+\,\vert\, b\in K\}.$
Now, if $\Theta_g=\Theta_h$, then $h=g^+_b$, for some $b$. 
But then $h=f\wedge g$ with $f(x)=x-b$, that is, the action 
of $G^+$ can be represented in terms of the wedge operator.

Let us denote the discriminant of $f$ by $\Delta(f)$. 
In the following lemma we collect miscellaneous results
on the action of $G$, to be used in later sections.
From part $(iii)$ we find that the discriminant is a 
block invariant.

\begin{lemma} \label{lemma:ActionOfG}
Let $f$ and $g$ belong to a stable set, let  $n={\rm deg\/} f$
and let $\sigma_{a,b}$ be an arbitrary element of $G$.
The following holds:
\begin{enumerate}
\item[$(i)$] $\sigma_{a,b}(\Theta_f)=\Theta_{\sigma_{a,b}(f)}.$
\item[$(ii)$] $\sigma_{a,b}(f)\wedge \sigma_{a,b}(g) \,=\, 
     \sigma_{a^{n},0}\left(f\wedge g\right).$
\item[$(iii)$] $\Delta(\sigma_{a,b}(f))\,=\,\Delta(f)/a^{n\,(n-1)}$.
\item[$(iv)$] If $f,g\in\Theta$, then, for all $b\in K$, we have 
$f\wedge f^+_b=g\wedge g^+_b$.
\end{enumerate}
\end{lemma}

\proof
If $g\in\Theta_f$, then $\sigma^+(f)=g$ for some
$\sigma^+\in G^+$, whence
$$
\sigma(g)= \sigma\sigma^+(f)=(\sigma\sigma^+\sigma^{-1})\,\sigma (f).
$$
Because $G^+$ is a normal subgroup, we have that
$\sigma\sigma^+\,\sigma^{-1}\in G^+$, so that
$\sigma(g)$ is in the same block as $\sigma(f)$.
Furthermore, $\sigma$ is bijective, because (\ref{eq:ActionOfG}) 
defines a group action.
This proves the first assertion.
To prove $(ii)$, let $g(\alpha)=0$. With reference to (\ref{eq:Wedge})
and  (\ref{eq:ActionOfG}), we find
$$
\sigma_{a^n,0}(f\wedge g)(x)=
  \mathop{{\prod}'}_\alpha\,\left(x-a^{-n}f(\alpha)\right);
\hskip 30pt
(\sigma_{a,b}(g))(x)=\mathop{{\prod}'}_\alpha\,\left(x-a^{-1}(\alpha-b)\right).
$$
from which we obtain
$$
\left(\sigma_{a,b}(f)\wedge \sigma_{a,b}(g)\right)(x)\,=\, 
\mathop{{\prod}'}_\alpha\,\left[x-a^{-n}f\left(aa^{-1}(\alpha-b)+b\right)\right]
 \,=\, \mathop{{\prod}'}_\alpha\,(x-a^{-n}f(\alpha))
$$
as required.
Part $(iii)$ follows from the discriminant formula
for a monic polynomial \cite[chapter 1.4]{LidlNiederreiter}
$$
\Delta(f)=\prod_{1\leq i < j \leq n} (\alpha_i-\alpha_j)^2
$$
where $\alpha_k$ are the roots of $f$, and the fact
that $\sigma_{a,b}(f)$ has the same degree as $f$.

Now, let $g,f\in\Theta$, so that $g=f^+_c$, for some $c\in K$.
Then, from $(ii)$ above, we have
\begin{eqnarray*}
g\wedge g^+_b&=&\sigma_{1,c}(f)\wedge\sigma_{1,b}\sigma_{1,c}(f)\\
 &=& \sigma_{1,c}(f)\wedge\sigma_{1,c}\sigma_{1,b}(f)\\
 &=& f\wedge f^+_b.
\end{eqnarray*}
This proves $(iv)$.
\endproof

\subsection{Self-image of blocks}
\label{section:SelfImageOfBlocks}

The self-image $\W(\Theta)$ of a block $\Theta$,
given by equation (\ref{eq:SelfImage}), consists
of all polynomials of the type $f\wedge f^+_b$, 
with $f\in\Theta$ and $b$ a non-zero element of $K$. 
From lemma \ref{lemma:ActionOfG} $(iv)$, it follows that
$f\wedge f^+_b$ does not depend on the choice of $f$ 
in $\Theta$, and therefore the self-image of a block
is highly degenerate. In particular, we find sufficient 
conditions for degree invariance, and for the existence 
of wedge invariants.

\begin{proposition}\label{proposition:DegreeInvariance}
Let $\Theta$ be a block of degree $n$, where $n$ is prime
and char$(K)\not=n$. Then the self-image of $\Theta$ has degree $n$.
\end{proposition}

\proof
Let $\Theta=\Theta_f$, and let $\alpha$ be a root of $f$.
If $b\not=0$, the polynomial $f\wedge f^+_b$ has root 
$f(\alpha-b)=-n\,b\,\alpha^{n-1} + \cdots$, which is a
polynomial in $\alpha$ of degree $n-1$, since, by assumption,
$b\,n\not=0$. 
Because $n$ is prime, the degree of $f(\alpha-b)$ is either 
$n$ or 1.
In the latter case, if $f(\alpha-b)=d\in K$, then we have 
$f(\alpha-b)-d=0$, which is impossible since $\alpha$ is 
algebraic of degree $n$.
So the degree of $f\wedge f^+_b$ is the same as that of $f$.
\endproof

\begin{theorem} \label{theorem:RightLeftInvariants}
If a block $\Theta$ intersects its self-image, then $\Theta$
contains a polynomial which is a right-left invariant of 
two polynomials in the same block.
\end{theorem}

\proof
By assumption, a polynomial $f$ exists such 
that $f\wedge f_b^+=f_c^+$ for some $b,c\in K$,
with $b\not=0$.
Then, letting $\theta=f_c^+$, from lemma \ref{lemma:ActionOfG}$\,(iv)$ we obtain
\begin{equation}
\theta_{-b}^+\wedge \theta=\theta\qquad\quad
\theta\wedge\theta_b^+=\theta, 
\end{equation}
that is, $\theta$ is the desired right-left 
invariant of the wedge operator.
\endproof

If, with the above notation, one also has that
$\theta_{b}^+\wedge \theta=\theta$, then the self-intersection
of a block yields two stable 2-sets of type III,
namely $\{\theta,\theta_{\pm b}^+\}$.
One verifies that this is the case precisely when $\theta(\alpha+b)$
and $\theta(\alpha-b)$ are algebraic conjugates,
which always happens for quadratic polynomials
--- see section \ref{subsection:QuadraticGraphs}.

An instance of this phenomenon is shown in figure \ref{fig:Stable3Set}.
The quadratic extension graph of $\F_3$ consists of a single block,
which intersects (indeed, contains) its self-image, due to proposition
\ref{proposition:DegreeInvariance}, giving rise to the
right-left invariant $g$, and to two stable 2-sets of type III.

If we now let $\B$ be the set of blocks of an extension,
the mapping
\begin{equation} \label{eq:SelfMap}
\W_b: \B\to \B
\hskip 40pt 
\Theta_f\,\mapsto\,\Theta_{f\wedge f^+_b}
\hskip 40pt 
0\not=b\in K
\end{equation}
is well-defined.
%
By means of $\W_b$ we can construct an oriented graph, 
whose vertices are blocks, and where two blocks 
$\Theta$ and $\Theta'$ are joined by an
arc if $\Theta'$ belongs to the self-image of $\Theta$.
This graph, called the {\it block graph,} contains the
essential information on all self-interactions: it will 
be studied in the next section for the quadratic case.

\section{Quadratic polynomials}
\label{section:Quadratic}

Throughout this section, $f$ denotes a quadratic
monic irreducible polynomial over a field $K$,
with discriminant $\Delta(f)$.
We first describe the graphs of a quadratic extension 
(theorem \ref{theorem:Quadratic}), and characterize
the transitions to the ground field (proposition \ref{proposition:Transition}).

Then we construct stable two sets, and derive a three-dimensional 
skew-map describing the dynamics of 2-sets, with which we compute
periodic points.
It turns out that for periodic 2-sets to exist, the ground
field $K$ must contain certain roots of unity.
In this respect, the situation is not dissimilar from that 
of periodic orbits of the monomial maps $z\mapsto z^k$,
for which the $\Phi$-polynomials ---cf.~equation \ref{eq:PhiPoly}---
are cyclotomic polynomials.

\subsection{Graphs}
\label{subsection:QuadraticGraphs}
The quadratic extension graphs have a particularly simple form. 
 
\begin{theorem} \label{theorem:Quadratic}
A quadratic extension graph consists of a single cluster, 
whose distinct blocks have distinct discriminants.
If char$(K)\not=2$, then the block-graph is a complete graph, 
and the set of mappings $\W_b$, defined in (\ref{eq:SelfMap}),
form a group of permutations of blocks, isomorphic to the 
multiplicative group of $K^2$. 
The isomorphism associates $b^2\in (K^*)^2$ to 
the permutation sending the block $\Theta$ to $\W_{b/2}(\Theta)$.
\end{theorem}

The completeness of the block graph means that the diameter
of the extension graph is at most two, combining polynomials 
of degree 1 and 2. 
This holds in a strong sense, namely for every quadratic
polynomials $f$ and $g$ there exist elements $b$ and $c$ in $K$
such that $(f_b^+\wedge f)_c^+=g$.

We need a lemma:

\begin{lemma}
\label{lemma:QuadraticIdentities}
For every quadratic polynomial $f$ and every $b\in K$, the following holds
\begin{enumerate}
\item[$(i)$] $(f\wedge f_b^+)(x)\,=\,(x-b^2)^2-b^2\,\Delta(f)$.
\item[$(ii)$] $\Delta(f\wedge f_b^+)\,=\,(2\,b)^2\Delta(f)$.
\item[$(iii)$] If char$(K)\not=2$, then the polynomials $f\wedge f_a^+$ and 
$f_b^+\wedge f_{b+c}^+$ belong to the same block if and only if $a=\pm c$, 
in which case they coincide.
\end{enumerate}
\end{lemma}

\proof
The first two identities are verified by direct calculation.
To prove $(iii)$, we first reduce it to the case $b=0$, 
from lemma \ref{lemma:ActionOfG}$\,(ii)$. 
Then from $(i)$ above, we have that if $f\wedge f_a^+$ and 
$f_b^+\wedge f_{b+c}^+$ belong to the same block, then
for some $d\in K$ 
$$
(x+d-a^2)^2-a^2\,\Delta(f)\,=\,(x-c^2)^2-c^2\,\Delta(f)
$$
which implies that $d=0$ and $a^2=c^2$.
Conversely, if $a=\pm c$, from part $(i)$ above and lemma 
\ref{lemma:ActionOfG}$\,(ii)$,
we have that $f\wedge f_a^+=f_b^+\wedge f_{b\pm a}^+=
f_b^+\wedge f_{b+c}^+$.
\endproof
 
\bigskip
\noindent {\it Proof of theorem \ref{theorem:Quadratic}} \quad
If $\alpha$ is a root of $f$ and $g$ is an irreducible 
polynomial in the same extension, then the roots of $g$ are 
linear expressions in $\alpha$ with coefficients in $K$.
Thus, for some $a,b\in K$, we have that $\sigma_{a,b}(f)=g$,
i.e., there is a single cluster.
Now, from lemma \ref{lemma:ActionOfG}$\,(iii)$ we have
$\Delta(\sigma_{a,b}(f))\,=\,\Delta(f)/a^2$, so distinct
values of $a^2$ correspond to distinct discriminants.
Now let $S$ be the sum of the roots of $f(x)$.
One verifies that
$$
\sigma_{a,b}(f)=\sigma_{1,c}(\sigma_{-a,b}(f)) 
\hskip 30pt c={2b-S\over a}
$$
which shows that $\sigma_{a,b}(f)$ and $\sigma_{-a,b}(f)$ 
belong to the same block, and so blocks are parametrized by
discriminants.

Let char$(K)\not=2$.
To show that the block graph is complete, we consider the equation
$$
\W_b(\Theta_f)=\Theta_{f\wedge f^+_b}=\Theta_g
$$
where $f$ and $g$ are given quadratic irreducible polynomials,
with $E(f)=E(g)$. We look at discriminants.
Because our extension is separable, and all discriminants 
are square multiples of the field discriminant, we have
$\Delta(g)=k^2\Delta(f)$, for some $0\not=k\in K$.
This, together with lemma \ref{lemma:QuadraticIdentities}$(ii)$,
gives the equation $k^2=(2b)^2$, which can be solved for $b$, since
char$(K)\not=2$. Matching discriminants suffices, since we have
seen that discriminants identify blocks. 

We show that $\W_b$ is injective.
If $\Theta$ and $\Theta'$ are distinct blocks, they 
have distinct discriminants, $\Delta$ and $\Delta'$, say. 
But then
$$
\Delta(\W_b(\Theta_f))=2b^2\Delta\not=2b^2\Delta'=\Delta(\W_b(\Theta_g)) 
$$
which shows that $\W_b(\Theta)\not=\W_b(\Theta')$, as desired.
To prove surjectivity, we must solve for $f$ the equation 
$\W_b(\Theta_f)=\Theta_g$, for given $b$ and $g$. 
From what was proved above, this amounts to find $f$ 
such that $\Delta(f)(2b)^2=\Delta(g)$;
since char$(K)\not=2$, we can take $f=\sigma_{2b,0}(g)$.

Thus, if char$(K)\not=2$, each non-zero value of $b$ defines
a permutation of blocks. 
Now consider the mapping 
$$
\mu:\,b^2\mapsto \W_{b/2}\hskip 40pt 0\not= b\in K
$$
sending $(K^*)^2$ to the symmetric group on $\B$.
The choice of $b$ among the square roots of $b^2$ is
irrelevant, due to lemma \ref{lemma:QuadraticIdentities}$(iii)$.
From the same lemma, part $(ii)$, we see that 
$\mu$ associates to $b^2$ the permutation sending the block 
of discriminant $\Delta$ to that of discriminant $b^2\Delta$.
Keeping this in mind, we find that
$$
\mu((bc)^2)=\W_{bc/2}=\W_{b/2}\circ\W_{c/2}=\mu(b^2)\mu(c^2)
$$
e.g., $\mu$ is a group homomorphism. Its kernel is trivial,
$\mu(1^2)=\W_{1/2}$, and so $\mu$ defines a faithful action 
of $(K^*)^2$ on $\B$.
\endproof

Thanks to theorem \ref{theorem:RightLeftInvariants},
the identity permutation $\W_{1/2}$ maps the
whole block to a right-left invariant of the wedge 
operator, which, from lemma \ref{lemma:QuadraticIdentities}$(i)$,
is given by
\begin{equation}\label{eq:Centre}
\theta(x)=\left(x- {1\over 4}\right)^2 -{\Delta\over 4}
\end{equation}
where $\Delta$ is the block discriminant. 
From lemma \ref{lemma:QuadraticIdentities}$(iii)$, this
polynomial is unique: we call it the {\it centre\/} of the block.
Because $\theta=\theta\wedge\theta^+_{1/2}$ $=\theta\wedge\theta^+_{-1/2}=
\theta^+_{1/2}\wedge\theta$ $=\theta^+_{-1/2}\wedge\theta$,
we obtain the following type III stable sets
\begin{equation}\label{eq:Stable2Sets}
\{\theta,\theta^+_{1/2}\}
\hskip 30pt
\{\theta,\theta^+_{-1/2}\}
\hskip 30pt {\rm char}(K)\not=2
\end{equation}
which are distinct and in bi-unique correspondence with the
block discriminants $\Delta\in K$.
An example is given in figure \ref{fig:Stable3Set},
with $g=\theta$, $h=\theta_{1/2}$, and $f=\theta_{-1/2}$.
(Note however, that the set $\{\theta,\theta^+_{1/2},\theta^+_{-1/2}\}$,
which is stable in figure \ref{fig:Stable3Set},
is not stable in general.)
%

In the following table, we display all parametrized 
families of quadratic stable 2-sets, for char$(K)\not=2$.
The absence of quadratic stable 2-sets of type II results 
from theorem \ref{theorem:NoQuadraticStableSetsOfTypeII} below.

\bigskip
\goodbreak
\centerline{T{\footnotesize ABLE} I:\/
{\footnotesize QUADRATIC STABLE} 2-{\footnotesize SETS}}
\nobreak\bigskip\nobreak
\hfil\vtop{\baselineskip 30pt\halign{
\quad \hfil  $#$ \hfil &
\quad \hfil  $#$ \hfil &
\quad \hfil   #  \hfil &
\qquad\hfil  $#$ \hfil &
\quad \hfil  $#$ \hfil \quad \cr
f& g & type & \cr
\noalign{\vskip 5pt\hrule\vskip 5pt}
x^2+r & x^2-x+r & I & f=\Phi_{2,g} & g=\Phi_{1,f}\cr
x^2+{1\over 2}x+r &
   x^2-{1\over 2}x+r & III & f|g^2 & g=\Phi_{1,f}\cr
x^2-{3\over 2}x+{1\over 2}+r &
   x^2-{1\over 2}x+r & III & f|g^2 & g=\Phi_{2,f}\cr
}}\hfil

\bigskip
All polynomials are irreducible over $\Q(r)$, where $r$ is 
regarded as an indeterminate; if instead $r$ is as a specific 
element of $K$, then irreducibility must be checked.
The last two columns describe the mutual relation between 
$f$ and $g$, in the notation of section \ref{section:Invariance}.
The type I set is of the form (\ref{eq:TypeI}), with one
of the two polynomials having zero middle coefficient, as 
explained in section \ref{subsection:Periodic2Sets}.
The two type III sets are rooted at the block centre
$g=\theta$, and correspond to the two possible 
permutations of the roots of $g$ by $f$
(cf.~equations (\ref{eq:Centre},\ref{eq:Stable2Sets}), 
with $r=(1-\Delta)/4$). 

The last result of this section characterizes some 
transitions to the ground field.
For char $K\not=2$, this characterization is complete.

\begin{proposition}\label{proposition:Transition}
For every quadratic polynomials $f$ and $g$, and every 
$b\in K$, if $E(f)$ has more than one block, then the 
polynomial $(f\wedge f_b^+)\wedge (g\wedge g_b^+)$ 
has degree one.
Conversely, if char$(K)\not=2$ and $h\wedge l$ has degree 1, 
with $h$ and $l$ quadratic not belonging to the same block, 
then there exists $c\in K$ such that
$$
h\,=\,(\theta_h)^+_c \qquad\quad
l\,=\,(\theta_l)^+_c
$$
where $\theta_h$ and $\theta_l$ are the centres of the 
respective blocks.
\end{proposition}

Under the above assumptions, we have $\W^2(E(f))\not\subset E(f)$,
that is, some polynomials in the second self-image of a quadratic 
extension collapse onto the ground field.
 
\medskip
\proof
If $\alpha$ is a root of $l=g\,\wedge\, g^+_b$, then, 
from lemma \ref{lemma:QuadraticIdentities}$\,(i)$,
we have that $(\alpha-b^2)^2=b^2\,\Delta(g)$, 
and therefore, letting $h=f\wedge f_b^+$, the root $h(\alpha)$ 
of $h\wedge l$ is given by
\begin{equation} \label{eq:RootOfTransition}
h(\alpha)\,=\,b^2\,(\Delta(g)-\Delta(f))\, \in\, K.
\end{equation}
This proves the first statement.
Conversely, assume that $h\wedge l$ has degree 1, and
let $\theta_h=h\wedge h^+_{1/2}$ and $\theta_l=l\wedge l^+_{1/2}$ 
be the centres of the respective blocks. 
Then $\theta_h\wedge \theta_l$ has degree 1, from the above.
Solving for $c$ the equation $h=(\theta_h)^+_c$, 
gives $c=B/2+1/4$, where $B$ is the middle coefficient of $h(x)$.
If $r\in K$ is the root of $h \wedge l$, 
then $h(x)-l(x)=r$, and hence $\Delta(h)-\Delta(l)=-4\,r$.
The equation $l=(\theta_l)^+_{d}$ now reads
$$
    \left(x+d-{1\over 4}\right)^2 -{\Delta(l)\over 4}=
    \left(x+d-{1\over 4}\right)^2 -{\Delta(h)\over 4} - r
    =l(x) = h(x) -r
$$
with the solution $d=B/2+1/4=c$.
\endproof

\subsection{Periodic 2-sets}
\label{subsection:Periodic2Sets}

We describe the periodic behaviour of 2-sets $\{f,g\}$, 
where $f$ and $g$ are quadratic polynomials of discriminants 
$\Delta(f)$ and $\Delta(g)$, and ${\rm char}(K)\not=2$. 
Let $\W(\{f,g\})=\{f',g'\}$, with $f'=f \wedge g$ and $g'=g \wedge f$. 
Since $f'(0)=g'(0)=\hbox{Res}(f,g)$ (see remark following
equation (\ref{eq:Wedge})), without loss of generality, we let
$$
f=x^2+bx+r \hskip 40pt g=x^2+cx+r.
$$
The corresponding primed coefficients are computed as
$$
b'=c(b-c)
\hskip 30pt
c'=-b(b-c)
\hskip 30pt
r'=r(b-c)^2.
$$
Defining
\begin{equation}\label{eq:uv}
u=b-c
\hskip 40pt
v=b+c
\end{equation}
we obtain a three-dimensional skew map over $K$
\begin{equation}\label{eq:AuxiliaryMap}
\Psi:K^3\rightarrow K^3 \hskip 40pt
(u,v,r)\,\mapsto\,(\pm uv,-u^2,ru^2)
\end{equation}
where the change of sign corresponds to exchanging $f$ and $g$, 
from (\ref{eq:uv}). 
Iteration gives (ignoring sign change)
\begin{equation}\label{eq:Psi}
\Psi^{t}(u,v,r)\,=\,\cases{
(z_t/u,-z_t/v,r(z_t/uv)^2) & $t$ odd\cr 
\noalign{\vskip 1pt}
(-u z_t,-vz_t,rz_t^2) & $t$ even\cr 
}\hskip 30pt t>0
\end{equation}
where 
\begin{equation}\label{eq:ze}
z_t=(u^2v)^{e_t}
\hskip 30pt
e_t={1\over 3}\left[2^t+(-1)^{t+1}\right]
\hskip 30pt
t=1,2,\ldots
\end{equation}
The sequence $e_t=1,1,3,5,11,21,43,85,171,341,\ldots$ satisfies 
the recursion relation $e_{t+1}=e_t+2e_{t-1}$,
with initial conditions $e_1=e_2=1$.

From (\ref{eq:AuxiliaryMap}) we find
$$
(u,0,r)\,\mapsto\, (0,-u^2,ru^2) \,\mapsto\, (0,0,0)
$$
and so, in order not to collapse to the trivial solution,
we must have $uvr\not=0$.
Since the discriminants of $f$ and $g$ evolve as
\begin{equation}\label{eq:Discriminants}
\Delta(f')=u^2\Delta(g)
\hskip 40pt
\Delta(g')=u^2\Delta(f)
\end{equation}
one sees that if $f$ and $g$ are irreducible (in particular, $r\neq0$),
and their discriminants are distinct ($uv\neq0$), these properties
are preserved along the orbit.

The periodic points equation reads
\begin{equation}\label{eq:PeriodicPoints}
\Psi^{t}(u,v,r)\,=\,(\pm u,v,r).
\hskip 30pt t>0
\end{equation}
We have two cases, depending on the parity of $t$.
To unify the notation, we define
\begin{equation}\label{eq:dt}
d_t\,=\,\cases{
3e_t-2=2^t-1  & $t$ odd\cr 
\noalign{\vskip 1pt}
e_t=(2^t-1)/3  & $t$ even\cr 
}\hskip 30pt t>0,
\end{equation}
with $e_t$ as in (\ref{eq:ze}). 

$i)$ Odd period. When $t$ is odd, using (\ref{eq:Psi}), we find no 
irreducible solution corresponding to the positive 
sign in (\ref{eq:PeriodicPoints}).
In particular, for $t=1$, we have
\begin{theorem}\label{theorem:NoQuadraticStableSetsOfTypeII}
If Char$(K)\neq2$, there are no quadratic stable sets of type II over $K$.
\end{theorem}
For the negative sign in (\ref{eq:PeriodicPoints}), we find 
\begin{equation}\label{eq:OddPeriod}
v^{d_t}+1=0
\hskip 30pt
u=\pm v
\end{equation}
for one choice of sign, as changing sign amounts to exchanging polynomials. 
The sign alternates along each orbit, and in every element of a cycle 
of odd period, precisely one of the two polynomials has zero middle coefficient.
We see that $v=\zeta$, where $\zeta$ is a $2d_t$-th root 
of unity (for background references on roots of unity,
see, e.g., \cite[page 39]{McCarthy} or \cite[chapter 27]{Hasse}).
For every root of unity, there is a one-parameter family of solutions,
parametrized by $K\setminus K^2$, which corresponds to varying $r$
while keeping the polynomials irreducible.
When $t=1$, we recover the type I stable set displayed in Table I.
Because $d_t>1$ for $t>2$, odd cycles with period greater than 
one can exist only if $K$ contains non-trivial roots of unity.

$ii)$ Even period. When $t$ is even, we find no solution of 
(\ref{eq:PeriodicPoints}) corresponding to the negative sign 
(unless char$(K)=2$, which we have excluded). 
For the positive sign, we obtain
\begin{equation}\label{eq:EvenPeriod}
(u^2v)^{d_t}+1=0
\end{equation}
that is, $u^2v=\zeta$, where $\zeta$ is a $2e_t$-th root of unity.
For each solution of this equation, we have again a one-parameter
family of periodic points.
Solving (\ref{eq:EvenPeriod}) for $t=2$, gives the 2-cycle
$\{f_0,g_0\}\longleftrightarrow\{f_1,g_1\}$, 
irreducible over $\Q(r,s)$
\bigskip

\hfil\vtop{\baselineskip 35pt\halign{
$\displaystyle #$ \hfil & \qquad
$\displaystyle #$ \hfil\cr
 f_0=x^2-{s^3+1\over 2\,s}\,x+r &
 g_0=x^2-{s^3-1\over 2\,s}\,x+r\cr
 f_1=x^2+{s^3-1\over 2\,s^2}\,x+{r\over s^2} &
 g_1=x^2-{s^3+1\over 2\,s^2}\,x+{r\over s^2}.\cr
}}\hfil

\medskip\noindent
As above, even cycles with period greater than 
two require non-trivial roots of unity.

\bigskip

In the above construction, if $\zeta$ is a primitive 
$2d_t$-th root of unity, then the period $t$ is minimal.
This follows from the fact that if $t'$ is a proper divisor of $t$, 
then $d_{t'}$ is a proper divisor of $d_t$ ---see equation (\ref{eq:dt}). 
So 2$d_t$-th roots of unity are needed to build all cycles of 
minimal period $t$. However, orbits of the same minimal period may 
originate from some roots of unity of lower order, and indeed the
problem of determining the minimal order for a given period 
---the $t$-cycles of minimal complexity--- is also of interest 
(cf.~\cite{Vivaldi:06}).
To compute such orders, we consider all divisors of $d_t$,
and remove from them the divisors of $d_{t'}$, for all 
$t'<t$ such that $t'$ divides $t$ and has the same parity as $t$.
Call $D_t$ the resulting set of divisors.
The parity condition is justified as follows.
If $t$ is even and $t'$ is odd (this is the only possibility), and 
$d\in D_t$ is a common divisor of $d_t$ and $d_{t'}$, then the set
of solutions of equation (\ref{eq:EvenPeriod}) for the exponent 
$d$ is larger than that of equation (\ref{eq:OddPeriod}) for 
the same exponent.

In the table below, we display the (half)-orders $d_t$ of the roots of 
unity needed to construct periodic 2-sets of quadratic polynomials,
for all periods $t\leq 14$.

\goodbreak
\centerline{T{\footnotesize ABLE} II: \/
{\footnotesize ROOTS OF UNITY FOR QUADRATIC PERIODIC 2-SETS}}

\nobreak\vspace*{10pt}
\hfil\vtop{\baselineskip 13pt\halign{
\hskip 5pt
\hfil $#$ \hfil &\qquad
      $#$ \hfil \hskip 5pt\cr
\mbox{period } t & \mbox{order } d_t\cr
\noalign{\vskip 5pt\hrule\vskip 5pt}
 1 & 1 \cr
 2 & 1 \cr
 3 & 7 \cr
 4 & 5 \cr
 5 & 31 \cr
 6 & 21,7,3 \cr
 7 & 127 \cr
 8 & 85,17 \cr
 9 & 511,73 \cr
10 & 341,31,11 \cr
11 & 2047,89,23 \cr
12 & 1365,455,273,195,105,91,65,39,35,15,13 \cr
13 & 8191 \cr
14 & 5461,127,43\cr
}}\hfil

\bigskip
Of note are the large fluctuations of arithmetical origin,
and the close relation between the fields $K$ needed to construct
these periodic sets, and the fields that contain the periodic points
of the monomial maps $z\mapsto z^2$, which are $(2^t-1)$th roots of unity.

\section{Finite fields}
\label{section:FiniteFields}

In this section we consider self-interacting polynomials
over a finite field $K=\F_q$ with $q$ elements, where 
$q=p^k$, $p$ prime (for background information on finite 
fields, see \cite{LidlNiederreiter}).
Because a finite field has a unique extension of 
degree $n$ for any $n$, any stable set consisting 
of polynomials of bounded degree is also finite. 
Here we address some natural counting questions 
(number and size of blocks, number of stable 2-sets, 
number of periodic orbits, etc.).
Furthermore, we construct explicitly the periodic quadratic 
2-sets described in the previous section, and investigate
numerically the occurrence of stable sets of higher degree.

We denote the stable set of all irreducible polynomials 
of degree $n$ over $\F_q$ by $E(q^n)$, without reference 
to polynomials ---cf.~equation (\ref{eq:E}). 
Clearly, each block of $E(q^n)$ contains at most $q$ 
polynomials, but it may have fewer of them, and in some cases 
a block may even consist of a single polynomial, a so-called 
affine $q$-polynomial \cite[chapter 3.4]{LidlNiederreiter}.
However, if $n$ is coprime to $q$, the block size
is maximal, and we have

\begin{theorem} \label{theorem:NumberOfBlocks} 
If gcd$(q,n)=1$, then the number of blocks in the extension 
graph of degree $n>1$ over $\Fq$ is given by
\begin{equation}\label{eq:NumberOfBlocks}
{1\over nq}\,\sum_{d|n}\,\mu(d)\,q^{n/d}
\end{equation}
where $\mu$ is the M\"obius function.
\end{theorem}

The proof of this theorem will require the following lemma.
\begin{lemma} \label{lemma:FixedPoints}
Let $g$ be an irreducible polynomial over a field $K$.
If $g=g^+_b$, for some $b\not=0$, then char$(K)>0$ and
the degree of $g$ is divisible by char$(K)$.
\end{lemma}

\proof 
We first show that if $\alpha$, $\alpha+b$ and $\beta$ are
roots of $g$, so is $\beta+b$. Let $\tau$ be an element
of the Galois group of $g$, sending $\alpha$ to $\beta$.
Then $\beta+b=\tau(\alpha)+b=\tau(\alpha+b)$, showing that
$\beta+b$ is conjugate to $\alpha+b$.
Now, let $H$ be the collection of elements $b$ of 
$K$ for which $\alpha+b$ is a root of $g$. 
Then $H\not=\{0\}$, by hypothesis, and it is an additive group,
as seen from repeated applications of the above argument.
It follows that char$(K)>0$ 
and that the subgroup $\langle\,b\rangle$ of $H$ has order 
$p:={\rm char}(K)$, and so the order of $H$ is divisible by $p$.
If $\beta$ is another root of $g$ not of the form $\alpha+b,\,
b\in K$, then the corresponding group has the same order as $H$, 
again from the above argument.
Repeating this procedure until all roots of $g$ are accounted 
for, yields the result.
\endproof 

\noindent {\it Proof of theorem \ref{theorem:NumberOfBlocks}}.
The number of irreducible polynomials of degree $n$ over the finite
field $\F_q$ is given by \cite[theorem 3.25]{LidlNiederreiter}
\begin{equation}\label{eq:NumberOfIrreduciblePolynomials}
\#E(q^n)={1\over n}\,\sum_{d|n}\,\mu(d)\,q^{n/d}.
\end{equation}
From lemma \ref{lemma:FixedPoints}, if $q$ and $n$ are coprime, no 
irreducible polynomial of degree $n$ over $\Fq$ can have roots 
differing by elements of $\Fq$. So if $f\in E(q^n)$, the $q$ 
polynomials $f_a^+, \, a\in\Fq$ are all distinct, and form a block of order $q$.
\endproof

Our next task is to count quadratic stable sets.
The following lemma will be needed to check the simultaneous
irreducibility of pairs of quadratic polynomials.

\begin{lemma} \label{lemma:NN} 
Let $q$ be odd, and let $a$ be a non-zero element of $\F_q$.
If $q\equiv 3\mod{4}$, then there are $(q-3)/4$ values of $x\in\F_q$ 
such that $x$ and $x+a$ are both non-zero and both non-squares. 
If $q\equiv 1\mod{4}$ then the number of such values of $x$ is
equal to $(q-1)/4$ if $a$ is a square, and to $(q-5)/4$ if $a$ 
is not a square.
\end{lemma}

\proof The case $q=3$ is trivial, so we assume $q>3$.
Consider the polynomial $L(x)=F(x)\,G(x)\,H(x)$ where
$$
F(x)=x^{(q-1)/2}-(x+a)^{(q-1)/2};  \hskip 10pt
G(x)=x^{(q-1)/2}+(x+a)^{(q-1)/2};  \hskip 10pt
H(x)=x(x+a).
$$
From Euler's criterion, for any $x\in\F_q$, precisely
one of $F(x)$, $G(x)$, $H(x)$ is zero: if $H(x)$ is non-zero,
then $F(x)$ is zero when $x$ and $x+a$ have the same quadratic
character (they are both squares or non-squares) and 
non-zero otherwise, and conversely for $G(x)$.
Furthermore, deg$F=(q-3)/2$, deg$G=(q-1)/2$, deg$H=2$,
and hence deg$L=q$, so that all roots of $L$ are distinct.
It follows that $L$ is a constant multiple of $x^q-x$.
So the $(q-3)/2$ roots of $F$ are the values of $x$ for which $x$
and $x+a$ are both non-zero, and have the same quadratic character.

Consider the involution $x\mapsto \iota(x)=-(x+a)$. We have two cases.

$i)$\/ If $q\equiv 3\mod{4}$, then $F(\iota(x))=F(x)$, and moreover $-1$ is not a square. 
Thus, if $\alpha$ is a root of $F$ so is $\iota(\alpha)$, and these two
roots have opposite quadratic character. For this reason the involution $\iota$
cannot fix any root of $F$, and hence exactly half of such values of $\alpha$,
$(q-3)/4$ in number, are such that both $\alpha$ and $\alpha+a$ are non-squares.

$ii)$\/ If $q\equiv 1\mod{4}$, then $G(\iota(x))=G(x)$, and $x$ and $-x$ have
the same quadratic character. If $\alpha$ is a root of $G$, then 
$\iota(\alpha)$ is another root of $G$ with opposite quadratic 
character. So $\iota$ cannot fix any root of $G$, and therefore half of the
roots of $G$ are squares, and half are non-squares. The remaining
$(q-1)/4$ non-squares are subdivided between the roots of $F$ and
the root $-a$ of $H$, whose quadratic character is the same as $a$. 
One sees that if $a$ is a square, then there are $(q-1)/4$
roots of $F$ which are non-squares, and if $a$ is not a square 
then the number of such roots is $(q-1)/4-1=(q-5)/4$.
\endproof

We can now count quadratic stable sets for odd $q$.

\begin{theorem}\label{theorem:NumberOfStable2Sets}
Let $q$ be odd, and let $N_i(q)$ be the number of stable quadratic 2-sets 
of type $i$ ($i=I,II,III$) over $\F_q$. We have
\begin{enumerate}
\item [] $N_I(q)=\cases{(q-1)/4& if $q\equiv 1\mod{4}$ \cr
                        (q-3)/4& if $q\equiv 3\mod{4}$;\cr}$
\item [] $N_{II}(q)=0$;
\item [] $N_{III}(q)=q-1$.
\end{enumerate}
\end{theorem}

\proof The case of type II is a specialization of theorem 
\ref{theorem:NoQuadraticStableSetsOfTypeII}.
For type III, theorem \ref{theorem:NumberOfBlocks}
for $q$ odd and $n=2$ gives gives $(q-1)/2$ blocks.
Now, every quadratic block has a centre, which gives rise to two 
distinct stable sets of type III according to equation (\ref{eq:Stable2Sets}). 

For type I, with reference to table I, we have to verify the simultaneous
irreducibility of $f$ and $g$, whose discriminant is $-4r$ and $1-4r$, respectively.
Because $r$ is an arbitrary non-zero element of $\F_q$ and $q$ is odd, each 
discriminant assumes $q-1$ distinct values. 
Our result now follows from lemma \ref{lemma:NN} with $a=1$.
\endproof

We finally turn to periodic 2-sets. 

\begin{theorem}\label{theorem:NumberOfPeriodic2Sets}
For odd $q$, the number of periodic 2-sets of degree 2 
over $\F_q$ is at most $(q-1)^2\cdot (q-3)/8$.
\end{theorem} 

\proof We use the results and notation of section 
\ref{subsection:Periodic2Sets}. A periodic 2-set $\{f,g\}$
of degree 2 has the form
\begin{equation}\label{eq:PeriodicSet}
f=x^2+bx+r,\hskip 20pt g=x^2+cx+r,
\hskip 20pt uvr=(b^2-c^2)r\neq 0.
\end{equation}
Interchanging the polynomials changes the sign of $u$,
and so the number of pairs $(u,v)$ to be considered
is equal to $(q-1)^2/2$.

We now apply lemma \ref{lemma:NN}, considering that the 
difference between the discriminants of $f$ ad $g$ is $uv$.
There are two cases.

$i)$\/ If $q\equiv 3\mod{4}$, then for every pair $(u,v)$ 
there are $(q-3)/4$ values of $r$ for which $f$ and $g$
are both irreducible. 

$ii)$\/ If $q\equiv 1\mod{4}$, then, given $(u,v)$, 
the number of values of $r$ with the stated property 
is equal to $(q-1)/4$ if $uv$ is a square, and to 
$(q-5)/4$ if it is not. 
So the total number of irreducible pairs of eventually
periodic polynomials of the form (\ref{eq:PeriodicSet}) 
is given by
$$
\frac{(q-1)^2}{4}\left(\frac{q-1}{4}+\frac{q-5}{4}\right)
  =\frac{(q-1)^2(q-3)}{8}.
$$
This gives the result.
\endproof

We close this section by describing the construction of periodic 
2-sets of degree two, over a finite field of odd characteristic.
For the sake of brevity, we consider only the case of 
even period $t$, which is the most interesting.
For these cycles to exist, the field $\F_q$ must contain 
the $d_t$-th roots of unity, with $d_t$ given by (\ref{eq:dt}). 
These roots of unity belong to all finite fields only for $t=2$,
as described in section \ref{subsection:Periodic2Sets};
for $t>2$, the fields are restricted by the $t$-dependent 
condition $q\equiv 1\mod {d_t}$, which follows from the fact 
that the multiplicative group of a finite field is is cyclic
\cite[theorem 2.8]{LidlNiederreiter}.
From Dirichlet's theorem on arithmetic progressions
\cite[chapter 10]{Cohn}, we obtain
at once infinitely many fields $\F_q$ supporting orbits of a given 
even period.
Each of these fields contains $\phi(2d_t)=\phi(d_t)$ 
primitive $2d_t$-th roots of unity $\zeta$ ($\phi$ 
is Euler's function \cite[chapter 2]{Apostol}), which are 
constructed from a primitive element $\eta\in\F_q$, by 
letting $\zeta=\eta^i$ for all $i$ such that 
$(q-1)/{\rm gcd}(i,q-1)=2d_t$. 

Now, fix a root of unity $\zeta$ of order $2d_t$.
The equation $u^2v=\zeta$ can be solved for $u$ precisely
when $\zeta/v$ is a square, that is, when $v$ and $\zeta$
are both squares or non-squares.
Accordingly, we let $v=v_i=\eta^i$, with $i$ of
the appropriate parity, to obtain $(q-1)/2$ distinct values $v_i$.
For each $v_i$, we obtain two distinct solutions 
$(u,v)=(\pm\sqrt{\zeta/v_i},v_i)$. However, changing the
sign of $u$ corresponds to interchanging $b$ and $c$, that is,
interchanging polynomials. So we get $(q-1)/2$ distinct unordered
pairs $(b,c)$; for each pair we form all triples $(b,c,r)$
with the property that the discriminants $b^2-4r$ and $c^2-4r$ 
are simultaneously non-squares.
Using again lemma \ref{lemma:NN}, we find that there are
$\phi(d_t)(q-1)(q-3)/8$ points of minimal period $t$ 
associated to primitive $2d_t$-th roots of unity.

For illustration, let $t=6$. From Table II we find that
constructing a full set of $6$-cycles requires the 42nd 
roots of unity. We also see that the minimal order for that
period is $d_6=3$, and the smallest finite field containing 
3rd (hence 6th) roots of unity is $\F_7$. The above formula 
gives at most $\phi(3)(7-1)(7-3)/8=6$ points of minimal period 6, 
so there is just one 6-cycle over $\F_7$, $\{f_t,g_t\}$, $t=0,\ldots,5$,
which is displayed below.
\bigskip

$$
\hfil\vcenter{\baselineskip 15pt\halign{
\hskip 10pt\hfil $\displaystyle #$ \hfil & \qquad
$\displaystyle #$ \hfil & \qquad
$\displaystyle #$ \hfil\cr
 t & \qquad f_t & \qquad g_t \cr
\noalign{\vskip 5pt\hrule\vskip 5pt}
0 & x^2  + 2 x + 3 &  x^2  +   x + 3\cr
1 & x^2  +   x + 3 &  x^2  + 5 x + 3\cr
2 & x^2  +   x + 6 &  x^2  + 4 x + 6\cr
3 & x^2  + 2 x + 5 &  x^2  + 3 x + 5\cr
4 & x^2  + 4 x + 5 &  x^2  + 2 x + 5\cr
5 & x^2  + 4 x + 6 &  x^2  + 6 x + 6\cr
}}
\hskip 50pt 
K=\F_{7}.
\hfil
$$
We note that the difference of the middle coefficients 
(the variable $u$ in (\ref{eq:uv}) runs through the
entire multiplicative group of $\F_7$.

\subsection{Some experiments} \label{subsection:Experiments}
We have explored with Maple the occurrence of stable 
sets over various finite fields.
In the following table we display the number of
stable 2-sets for all extensions $E(p^n)$ of a prime 
fields $\F_p$ containing fewer than 500 polynomials.
Here $n$ is the degree of the extension,
$\#E$ denotes the number of irreducible polynomials
of degree $n$, while I, II, III denote the
type of stable set.

\begin{figure}[b]
\hfil\epsfig{file=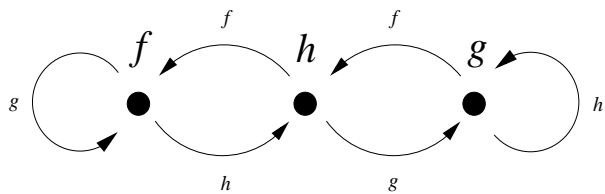,width=8cm,height=2.5cm}\hfil
\caption{\label{fig:Invariant3Set} 
The graph of the invariant 3-set $\{f,g,h\}$ of degree 4 over $\F_5$,
given in (\ref{eq:Invariant3Set}).
}
\end{figure}

\goodbreak
\centerline{T{\footnotesize ABLE} III: \/
{\footnotesize NUMBER OF STABLE} 2-{\footnotesize SETS}}
\nobreak\vspace*{10pt}
\hfil\vtop{\baselineskip 13pt\halign{
\hskip 5pt
\hfil $#$  &\qquad
\hfil $#$  &\qquad
\hfil $#$  &\qquad
\hfil $#$  &\qquad
\hfil $#$  &\qquad
\hfil $#$ \hskip 5pt\cr
 p & n & \#E &\rm I &\rm  II &\rm  III \cr
\noalign{\vskip 5pt\hrule\vskip 5pt}
 2 & 3&  2& -& -& 1 \cr
 2 & 6&  9& -& -& 1 \cr
 2 & 7& 18& -& -& 1 \cr
 2 & 9& 56& -& -& 2 \cr
 2 &10& 99& -& 1& - \cr
 2 &12&335& -& -& 1 \cr
 3 & 2&  3& -& -& 2 \cr
 3 & 3&  8& -& 3& - \cr
 3 & 4& 18& -& -& 5 \cr
 3 & 5& 48& 2& -& - \cr
 3 & 6&116& -& -&11 \cr
 3 & 7&312& -& 1& - \cr
 5 & 2& 10& 1& -& 4 \cr
 5 & 4&150& 2& 1&14 \cr
 7 & 2& 21& 1& -& 6 \cr
 7 & 3&112& 2& -& - \cr
11 & 2& 55& 2& -&10 \cr
11 & 3&440& -& 4& - \cr
13 & 2& 78& 3& -&12 \cr
17 & 2&136& 4& -&16 \cr
19 & 2&171& 4& -&18 \cr
23 & 2&253& 5& -&22 \cr
29 & 2&406& 7& -&28 \cr
}}\hfil

\bigskip\noindent
If a row in the table is missing (e.g., $p=2$, $n=2$), it means 
that there are no stable 2-sets in that extension. 
For $p\neq 2$, the data for quadratic extension follow from
theorem \ref{theorem:NumberOfStable2Sets}, which, in particular,
explains the absence of type II sets and the abundance of type III. 
Beyond the quadratic case, type II sets seem rare; however,
$n=3$ suffices, e.g., $f=x^3+x^2+2$, $g=x^3+2x^2+1$, over $\F_3$.
We observe that in the above table, for $p$ odd, type III stable 
sets occur only for extensions of even degree. 

Combinations of stable 2-sets may lead to interesting
invariant 3-sets, such as the one displayed in figure 
\ref{fig:Invariant3Set}.
It is given by
\begin{equation}\label{eq:Invariant3Set}
f=x^4+x+4,\quad g=x^4+2 x+4,\quad h=x^4+3 x+4
\hskip 30pt K=\F_5
\end{equation}
and it contains two stable 2-sets: $\{f,h\}$ (type II),
and $\{g,h\}$ (type III).

We close with some remarks and computations on periodic orbits.
The dynamics of self-interactions over a finite field is 
eventually periodic, and a natural problem is to determine 
the structure of periodic sets over a given field.
Even in its simplest form ---determining the period of 
the limit cycle of a quadratic 2-set--- this problem seems at 
least as difficult as that of computing the period of the squaring
map $x\mapsto x^2$ over a finite field, see comment at the end of 
section \ref{subsection:Periodic2Sets}. Some formulae and 
asymptotic expressions concerning the periodicity of the 
squaring map (and, more generally, of repeated exponentiation) 
have been obtained in \cite{ChouShparlinski}.

We are interested in a specific probabilistic phenomenon.
Given a set $z=\{f,g\}$ of two quadratic irreducible polynomials 
over $\Z$, with the same constant coefficient, and such that 
the sum and difference of their middle coefficients is non-zero, 
we consider the orbit of $z$ over the field $\F_p$, $p$ prime. 
This orbit is eventually periodic: how long are the transient 
and the period of the limit cycle? 

Here we are pursuing an analogy with Artin's problem on 
primitive roots \cite{RamMurty}, which we now describe 
from a dynamical systems perspective. Given an integer $a$,
not a square, we consider the map $x\mapsto ax\mod{p}$, 
for $p$ coprime to $a$. The period of the orbit of any 
non-zero point $x\in\F_p$ is given by $\rm{ord}_p(a)$,
the order of $a$ modulo $p$. We define the normalized
period $T(p)=\rm{ord}_p(a)/(p-1)$. Regarding $T$ as
a random variable, one is interested in its distribution
function ${\cal D}(x)$, which is the probability that
$T$ assumes a value not exceeding $x$ for a prime $p$
chosen at random. Such a probability is computed using 
the natural density over the primes.
Artin conjectured\footnote{The validity of this 
conjecture is now known to follow from the generalized 
Riemann hypothesis.}
that ${\cal D}$ exists, and is a step
function with positive steps at the reciprocal of each 
natural number: $1,1/2,1/3,\ldots$ (figure \ref{fig:Artin}).
Furthermore, ${\cal D}$ does not depend on $a$, as long as 
$a$ is square-free, and has only a mild $a$-dependence 
otherwise. 


\begin{figure}[h]
\vspace*{-10pt}
\hfil\epsfig{file=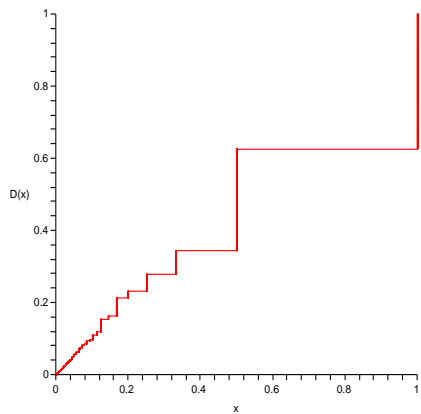,width=6cm,height=6cm}\hfil
\caption{\label{fig:Artin} 
Distribution function for the normalized period of the
orbits of the map $x\mapsto 2x \mod{p}$, computed over 
the first 20,000 odd primes.
The height of the step at $x=1$ is the so-called Artin's
constant: \/ $\prod_{p\geq 2}(1-(p^2-p)^{-1})=0.37395\ldots$ 
\cite[p.~303]{Ribenboim}.
}
\end{figure}

We have considered the orbit of the following 
pair of polynomials
\begin{equation}\label{eq:Artin}
z=\{f,g\}\hskip 40pt
f=x^2+x-1,\quad g=x^2+2x-1.
\end{equation}
The discriminants are $\Delta(f)=5,\Delta(g)=8$. Using 
quadratic reciprocity and the Chinese remainder theorem, 
we find that $f$ and $g$ are simultaneously irreducible 
modulo $p$ for $p\equiv{3,13,27,37}\mod{40}$. 
Excluding $p=3$ (for which the sum of the middle coefficients 
of $f$ and $g$ vanish), at all these primes, the orbit 
of $z$ consists of quadratic polynomials.
Accordingly we have considered, among the first 100,000 primes,
($p>3$) those belonging to the aforementioned residue classes
---25037 primes in all. 
For each prime $p$ we have computed the normalized period 
$T(p)=t(p)/(p-1)$, where $t(p)$ is the period of 
the limit cycle of the orbit of $z$ over $\F_p$. 

\begin{figure}[t]
\vskip -50 pt
\hfil\epsfig{file=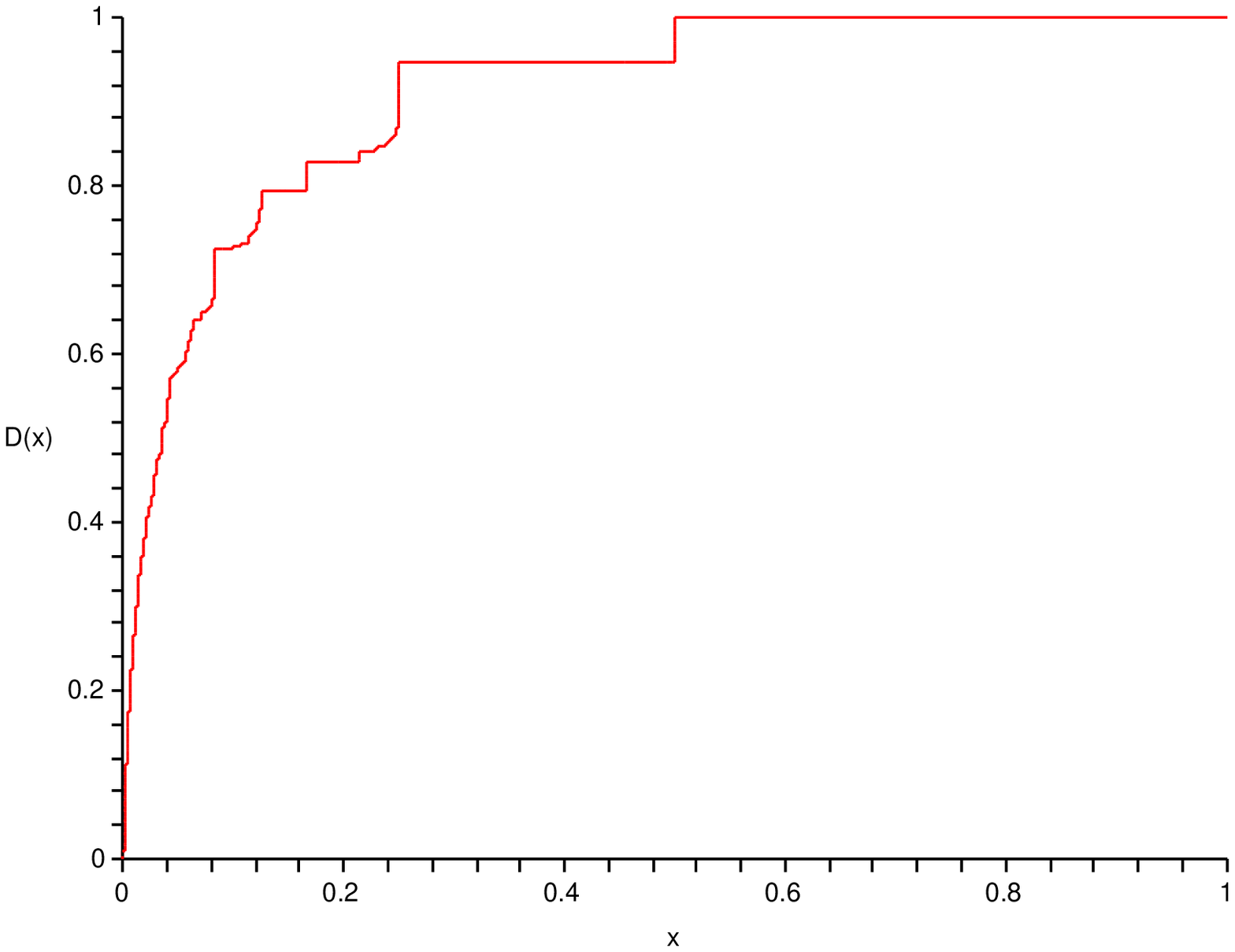,width=6cm,height=6cm}\hfil
\newline
\hfil\epsfig{file=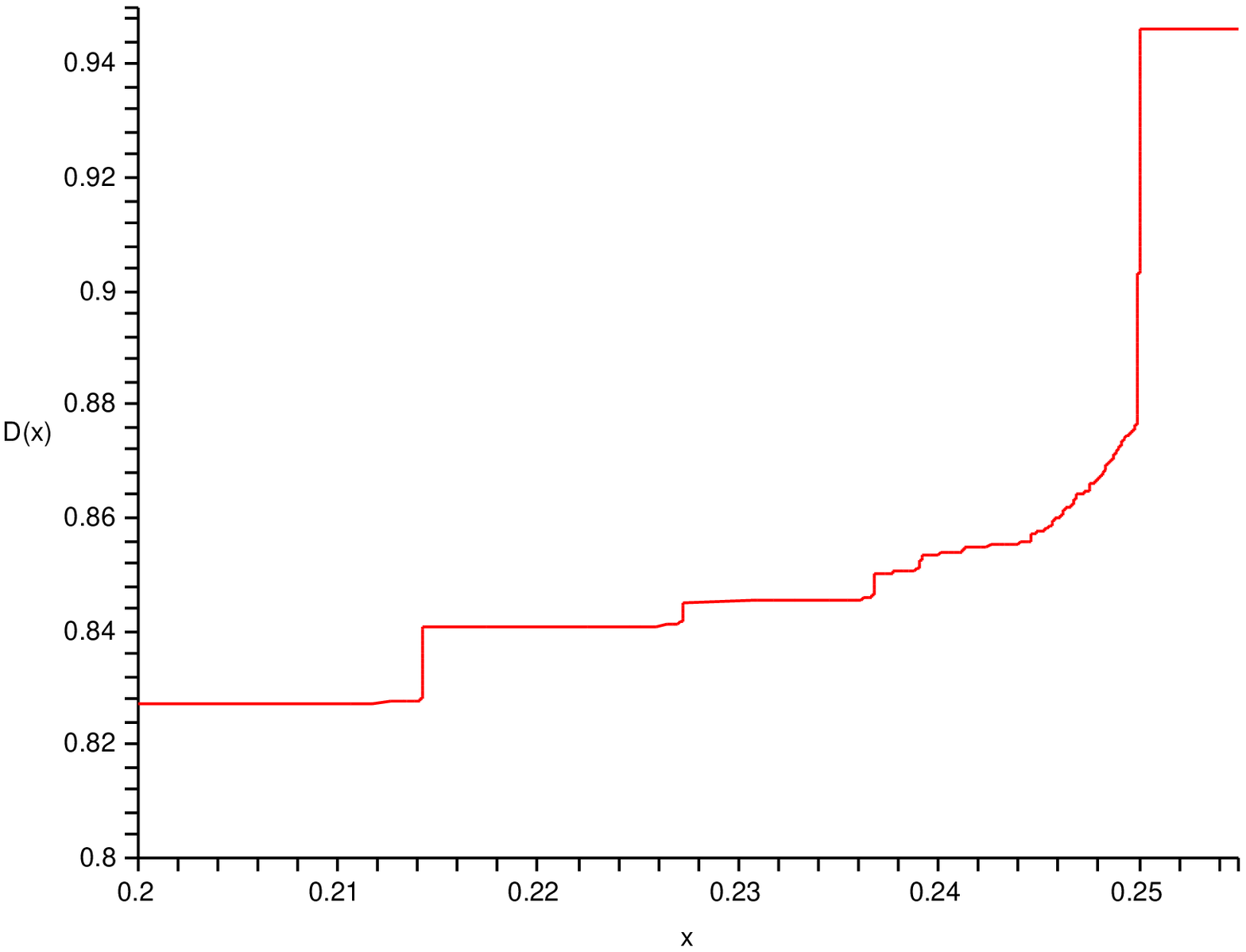,width=5cm,height=5cm}\hfil
\qquad
\hfil\epsfig{file=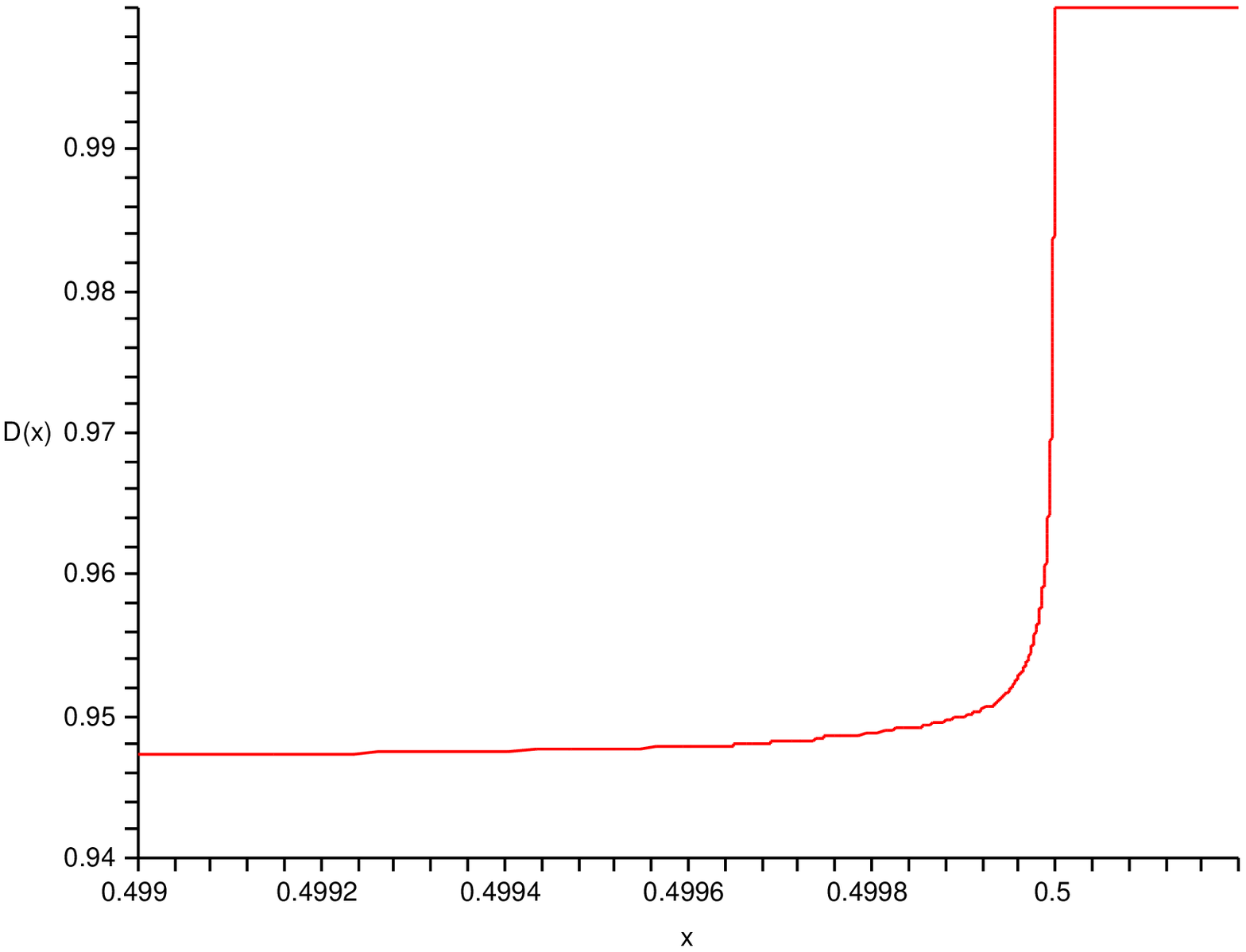,width=5cm,height=5cm}\hfil
\caption{\label{fig:Distribution} 
Top: distribution function for the normalized period of the
limit cycle of the orbit of the pair $z$ of equation (\ref{eq:Artin}).
The distribution function has been computed over a 
set of 25037 primes. 
Bottom: magnifications. The step at $x=1/4$ features 
secondary steps leading to it (left); that at $x=1/2$
has more regular climb (right).
}
\end{figure}

In figure \ref{fig:Distribution} we display the distribution 
function ${\cal D}(x)$. Its value reaches 1 at $x=1/2$, indicating 
that the cycle length does not exceed $(p-1)/2$ in a significant 
number of cases (in fact, we found $T(p)=1$ only for $p=163$).
This function has steps at the reciprocal of even integers.
Some steps have a clear sub-structure of secondary steps 
($x=1/4,1/8$, see figure \ref{fig:Distribution}, bottom left),
while others appear to be `smoother' 
($x=1/2,1/6$, see figure \ref{fig:Distribution}, bottom right).
There are no steps for odd denominators $x=1, 1/3, 1/5, 1/7$.
The orbit of $z$ was found to be either periodic 
(in very nearly $3/8$ of cases), or to have a transient 
of 1 or 2 (in $3/8$ and $1/4$ of cases, respectively).
This tight organization of transients was unexpected, and
we found it only for this specific value of $z$; in other 
examples we found instead a rapidly decaying distribution 
of transient lengths, consistent with the existence of 
a finite average transient. We remark that the finiteness of 
the average transient length has been proved for the squaring 
map (see \cite[theorem 2]{ChouShparlinski}).

We found that the distribution function does depend on the choice 
of initial conditions, although its basic structure remains the same.
The study of this function lies beyond the scope of this paper. 
Here we merely observe that, since $uvz\not=0$, 
the auxiliary map (\ref{eq:AuxiliaryMap}) can be transformed 
into an affine map of $(\Z/(q-1)\Z)^3$, using discrete logarithms.
In this setting, it should be possible to develop a qualitative 
analysis, although quantitative results are bound to be a lot 
more difficult.



\begin{thebibliography}{99}

\bibitem{Apostol}
 T. M. Apostol,
 Introduction to analytic number theory,
 Springer Verlag, New York
 (1976).

\bibitem{Bousch}
 T. Bousch,
 Sur quelques probl\`emes de la dynamique holomorphe,
 Ph.D.~thesis,
 Universit\'e de Paris-Sud, Centre d'Orsay
 (1992).

\bibitem{BatraMorton}
 A. Batra and P. Morton,
 Algebraic dynamics of polynomial maps on the algebraic closure of a finite field I,
 {\sl  Rocky Mountain J. of Math.}
 {\bf 24}
 (1994)
 453--481.

\bibitem{BatraMortonII}
 A. Batra and P. Morton,
 Algebraic dynamics of polynomial maps on the algebraic closure of a finite field II,
 {\sl  Rocky Mountain J. of Math.}
 {\bf 24}
 (1994)
 905--932.

\bibitem{ChouShparlinski}
 {W.-S. Chou and I. E. Shparlinski},
 {On thye cycle structure of repeated exponentiation modulo a prime}
 {\sl J. Number Theory}
 {\bf 107}
 {(2004)}
 {345--356}.

\bibitem{CohenHachenberger:99}
 {S. D. Cohen and D. Hachenberger},
 {Actions of linearized polynomials on the algebraic closure of a finite field},
 {\sl Contemporary Mathematics}
 {\bf 225}
 (1999)
 {17--32}.

\bibitem{CohenHachenberger:00}
 {S. D. Cohen and D. Hachenberger},
 {The dynamics of linearized polynomials},
 {\sl Proc. Edimb. Math. Soc.}
 {\bf 43}
 {(2000)}
 {113--128}.

\bibitem{Cohn}
 {H. Cohn},
 {\sl Advanced number theory},
 {Dover},
 {New York}
 {(1980)}.

\bibitem{Fontana:90}
  W. Fontana,
  Algorithmic Chemistry,
  in {\sl Artificial Life II,}
  SFI Studied in the Sciences of Complexity, vol X
  (1990)
  150--209.

\bibitem{FontanaBuss:94}
  W. Fontana and L. W. Buss,
  ``The arrival of the fittest'', toward a theory of biological organization,
  {\sl Bulletin of Mathematical Biology}
  {\bf 56}
  (1994)
  1--64.

\bibitem{Hasse}
 {H. Hasse},
 {\sl Number theory},
 {Springer-Verlag},
 {New York}
 {(2000)}.

\bibitem{KataokaKaneko:00}
 {N. Kataoka and K. Kaneko},
 {Functional Dynamics I: Articulation process,}
 {\sl Physica D}
 {\bf 138}
 (2000)
 225--250.

\bibitem{KataokaKaneko:01}
 {N. Kataoka and K. Kaneko},
 {Functional Dynamics II: Syntactic structure}
 {\sl Physica D}
 {\bf 149}
 (2001)
 174--196.

\bibitem{Koblitz:97}                                                            
 {N. Koblitz},                                                                  
 {\sl Algebraic aspects of cryptography},                                       
 {Springer-Verlag},                                                             
 {New York}                                                                    
 {(1997).}                                                                      

\bibitem{LidlNiederreiter}
 R. Lidl and H. Niederreiter,
 {\sl Finite fields},
 Encyclopedia of Math. and its Appl., vol. 20,
 Addison-Wesley, Reading, Mass.
 (1983).

\bibitem{Marcus}
 {D. A. Marcus},
 {\sl Number fields},
 {Springer-Verlag},
 {New York}
 {(1977)}.

\bibitem{McCarthy}
 {P. J. McCarthy},
 {Algebraic extensions of fields},
 {Dover},
 {New York}
 {(1966)}.

\bibitem{Morton:96}
  {P. Morton},
  {On certain algebraic curves related to polynomial maps},
  {\sl Compos. Math.}
  {\bf 103}
  {(1996)}
  {319--350}.

\bibitem{MortonPatel}
 P. Morton and P. Patel,
 The Galois theory of periodic points of polynomial maps,
 {\sl Proc. London Math. Soc.}
 {\bf 68}
 (1994)
 225--263.

\bibitem{MortonSilverman}
 P. Morton and J. H. Silverman,
 Periodic points, multiplicities, and dynamical units,
 {\sl J.~fuer reine und angew.~Math} 
 {\bf 461}
 (1995)
 {81--122}.

\bibitem{MortonVivaldi}
 P. Morton and F. Vivaldi,
 Bifurcations and discriminants for polynomials maps,
 {\sl Nonlinearity}
 {\bf 8}
 (1995)
 571--584.

\bibitem{Odoni:85}
 {R. W. K. Odoni},
 {The {G}alois theory of iterates and composites of polynomials},
 {\sl Proc. London Math. Soc.}
 {\bf 51}
 {(1985)}
 {385--414}.
 
\bibitem{RamMurty}
 {M. RamMurty},
 {Artin's conjecture for primitive roots},
 {\sl Math. Intelligencer}
 {\bf 10}
 {(1988)}
 {59--67}.

\bibitem{Ribenboim}
 {P. Ribenboim},
 {\sl The book of prime number records}
 Springer-Verlag,
 New York
 (1988).

\bibitem{Vivaldi}
 F. Vivaldi,
 Dynamics over irreducible polynomials,
 {\sl Nonlinearity}
 {\bf 5}
 (1992)
 941--960.

\bibitem{Vivaldi:06}
 F Vivaldi, 
 The arithmetic of discretized rotations,
 in {\sl $p$-adic Mathematical Physics,} 
 A. Y. Khrennikov, Z. Rakic, I. V. Volovich editors,
 AIP Conference Proceedings 826, AIP, 
 Melville, New York 
 (2006) 
 162--173.

\bibitem{VivaldiHatjispyros}
 F. Vivaldi and S. Hatjispyros,
 Galois theory of periodic orbits of rational maps,
 {\sl Nonlinearity}
 {\bf 5}
 (1992)
 961--978.

\bibitem{Waerden}
 {B. L. {van der Waerden}},
 {\sl Algebra},
 {Springer-Verlag},
 {New York}
 {(1991)}.

\end{thebibliography}
\end{document}